\documentclass[12pt,reqno]{amsart}

\usepackage{amssymb}

\usepackage{eucal}

\input{diagrams}

\font\sss=cmss8

\def\cG{{\mathcal G}}

\def\BZ{{\mathbb Z}}

\def\fg{{\mathfrak g}}

\def\fp{{\mathfrak p}}

\def\sC{\mbox{\sf C}}
\def\sD{\mbox{\sf D}}

\def\sK{\mbox{\sf K}}

\def\ast{{\textstyle *}}

\def\c{\operatorname{c}}

\def\D{\sD}

\def\Df{\sD^{\operatorname{f}}}

\def\fg{\operatorname{fg}}

\def\H{\operatorname{H}}

\def\Hom{\operatorname{Hom}}

\def\Mod{\mbox{\sf Mod}}

\def\opp{\operatorname{op}}
\def\pd{\operatorname{pd}}
\def\Proj{\mbox{\sf Pro}}

\def\RHom{\operatorname{RHom}}

\def\Z{\operatorname{Z}}

\numberwithin{equation}{part}

\newtheorem{Lemma}{Lemma}[section]
\newtheorem{Theorem}[Lemma]{Theorem}
\newtheorem{Proposition}[Lemma]{Proposition}

\theoremstyle{definition}

\newtheorem{Setup}[Lemma]{Setup}

\newtheorem{Construction}[Lemma]{Construction}
\newtheorem{Remark}[Lemma]{Remark}

\def\R{A}
\def\B{B}

\def\KProjR{\sK(\Proj\,\R)}
\def\KProjRopp{\sK(\Proj\,\R^{\opp})}
\def\KProjRsmall{\mbox{\sss K}{\scriptstyle (}\mbox{\sss Pro}\,{\scriptstyle \R)}}
\def\KModR{\sK(\Mod\,\R)}

\def\generatingset{\cG}
\def\generatingobj{G}

\def\fg{fi\-ni\-te\-ly ge\-ne\-ra\-ted}
\def\fp{fi\-ni\-te\-ly pre\-sen\-ted}

\def\cg{com\-pact\-ly ge\-ne\-ra\-ted}

\def\tr{tri\-an\-gu\-la\-ted}

\def\Rlm{$\R$-left-mo\-du\-le}
\def\Rrm{$\R$-right-mo\-du\-le}

\begin{document}

\title[Homotopy category]
{The homotopy category of complexes of projective modules}

\author{Peter J\o rgensen}
\address{Department of Pure Mathematics, University of Leeds,
Leeds LS2 9JT, United Kingdom}
\email{popjoerg@maths.leeds.ac.uk, www.maths.leeds.ac.uk/\~{ }popjoerg}




\begin{abstract} 
The homotopy category of complexes of projective left-modules over any
reasonably nice ring is proved to be a \cg\ \tr\ category, and
a duality is given between its subcategory of compact objects and the
finite derived category of right-modules.
\end{abstract}

\maketitle

\setcounter{section}{-1}
\section{Introduction}
\label{sec:introduction}

The last decade has seen \cg\ \tr\ categories rise to prominence.
Triangulated categories go back to Puppe and Verdier, but only later
developments have made it clear that the \cg\ ones are particularly
useful.  For instance, they allow the use of the Brown
Representability Theorem and the Thomason Localization Theorem, both
proved by Neeman in \cite{NeemanDuality}.  There are also results by
many other authors to support the case.

The standard examples of \cg\ \tr\ categories are the stable homotopy
category of spectra and the derived category of a ring.  Indeed, many
analogies between these two cases are captured by their common
structure of \cg\ \tr\ category, and this allows the transfer of
methods and ideas back and forth.

This paper adds to the collection of \cg\ \tr\ categories by showing
that if $\R$ is a reasonably nice ring, then the homotopy category of
complexes of projective \Rlm s, $\KProjR$, is \cg.

This may seem slightly surprising in view of \cite[app.\
E.3]{NeemanBook} which shows that the homotopy category of complexes
of all $\BZ$-modules, $\sK(\Mod\,\BZ)$, is not even well generated, a
weaker notion than \cg.  However, not only is $\KProjR$ compactly
generated; the subcategory of compact objects, $\KProjR^{\c}$, is very
nice, in that it is dual to the finite derived category of
\Rrm s, $\Df(\R^{\opp})$, whose objects are complexes with bounded
cohomology consisting of finitely presented modules.  My proofs of
these statements work when $\R$ is coherent and satisfies that each
flat \Rlm\ has finite projective dimension.

Most rings encountered in nature, such as noetherian rings, are
coherent.  The condition that each flat \Rlm\ has finite projective
dimension would appear less standard, but is in fact satisfied by
large classes of rings such as noetherian commutative rings of finite
Krull dimension (\cite[Seconde partie, cor.\ (3.2.7)]{RaynaudGruson}),
left-perfect rings (\cite[thm.\ P]{Bass}), and right-noetherian
algebras which admit a dualizing complex (\cite{PJfdpd}).  

The last of these cases includes many non-commutative algebras (see
\cite{WuZhangDualizing} and \cite{YekutieliZhang}), among them
noetherian complete semi-local PI algebras (\cite[cor.\
0.2]{WuZhangDualizing}) and filtered algebras whose associated graded
algebras are connected and noetherian and either PI, graded FBN, or
with enough normal elements (\cite[cor.\ 6.9]{YekutieliZhang}).

It is worth noting that if $\R$ has finite left and right global
dimension, then there is nothing new in my results.  In this case,
$\KProjR$ is equivalent to $\D(\R)$, the derived category of \Rlm s,
so $\KProjR$ is compactly generated.  Moreover, the subcategory of
compact objects $\KProjR^{\c}$ is equivalent to $\D(\R)^{\c}$, the
subcategory of compact objects of $\D(\R)$, and when $\R$ has finite
global dimension, $\D(\R)^{\c}$ is well known just to be the finite
derived category $\Df(\R)$, which is again dual to $\Df(\R^{\opp})$
under the functor $\RHom_{\R}(-,\R)$.

However, my results work for many rings which do not have finite
global dimension.

\section{Compact objects}
\label{sec:compact}

\begin{Setup}
\label{set:compact}
In this section, $\R$ is a right-coherent ring.
\end{Setup}

\begin{Construction}
\label{con:P}
Let $M$ be a \fp\ \Rlm.  This means that there is an exact sequence of
\Rlm s $Q_1 \rightarrow Q_0 \rightarrow M \rightarrow 0$ where
$Q_0$ and $Q_1$ are \fg\ projective \Rlm s.  

Hence there is an exact sequence of \Rrm s $0 \rightarrow M^{\ast}
\rightarrow Q_0^{\ast} \rightarrow Q_1^{\ast}$, where $(-)^{\ast}$
denotes the functor $\Hom(-,\R)$ which dualizes with respect to
$\R$.

Here $Q_0^{\ast}$ and $Q_1^{\ast}$ are \fg\ projective \Rrm s.  As
$M^{\ast}$ is the kernel of a homomorphism between them and as
$\R$ is right-coherent, it follows that $M^{\ast}$ is \fp.  Hence
$M^{\ast}$ has a projective resolution $P$ consisting of \fg\
projective \Rrm s.

Viewing $M^{\ast}$ as a complex concentrated in degree zero, there is
a canonical quasi-isomorphism 
\[
  P \stackrel{\pi}{\longrightarrow} M^{\ast}.  
\]
There is also a canonical homomorphism $M
\stackrel{\mu}{\longrightarrow} M^{\ast \ast}$ which I will view as
a chain map of complexes concentrated in degree zero, and so I can
consider 
\[
  M 
  \stackrel{\mu}{\longrightarrow} M^{\ast \ast}
  \stackrel{\pi^{\ast}}{\longrightarrow} P^{\ast}.
\]
\end{Construction}

\begin{Lemma}
\label{lem:Q1}
If $Q$ is a projective \Rlm, then
\[
  \begin{diagram}[labelstyle=\scriptstyle]
    \Hom_{\R}(P^{\ast},Q)
    & \rTo^{\Hom_{\R}(\pi^{\ast}\mu,Q)} &
    \Hom_{\R}(M,Q) \\
  \end{diagram}
\]
is a quasi-isomorphism.
\end{Lemma}

\begin{proof}
As $Q$ is projective, it is a direct summand in a free module, so it
is enough to prove the lemma when $Q$ is free.  But both $P^{\ast}$
and $M$ consist of \fp\ modules so when $Q$ is free,
and so has the form $\coprod \R$, then the coproduct can be moved
outside the $\Hom$'s, and so it is enough to prove the lemma for $Q =
\R$. 

There is a commutative diagram
\[
  \begin{diagram}[labelstyle=\scriptstyle]
    P & & \rTo^{\pi} & & M^{\ast} \\
    \dTo^p & & & & \vEq \\
    P^{\ast \ast} & \rTo_{\pi^{\ast \ast}} & M^{\ast \ast \ast} 
      & \rTo_{\mu^{\ast}} & M^{\ast} \\
  \end{diagram}
\]
where $p$ is the canonical chain map.  Since $P$ consists of \fg\
projective modules, $p$ is an isomorphism.  Also, $\pi$ is a
quasi-isomorphism by construction, so the diagram shows that the
composition $\mu^{\ast}\pi^{\ast \ast}$ is a quasi-isomorphism.

That is, the chain map
\[
  \mu^{\ast}\pi^{\ast \ast} 
  = (\pi^{\ast} \mu)^{\ast} 
  = \Hom_{\R}(\pi^{\ast} \mu,\R)
\]
is a quasi-isomorphism, and this proves the lemma in the case $Q = \R$
as desired.
\end{proof}

\begin{Lemma}
\label{lem:Q2}
If $Q$ is a complex of projective \Rlm s, then
\[
  \begin{diagram}[labelstyle=\scriptstyle]
    \Hom_{\R}(P^{\ast},Q)
    & \rTo^{\Hom_{\R}(\pi^{\ast}\mu,Q)} &
    \Hom_{\R}(M,Q) \\
  \end{diagram}
\]
is a quasi-isomorphism.
\end{Lemma}

\begin{proof}
The chain map $M \stackrel{\pi^{\ast}\mu}{\longrightarrow} P^{\ast}$
can be completed to a distinguished triangle
\[
    M \stackrel{\pi^{\ast}\mu}{\longrightarrow} 
    P^{\ast} \longrightarrow C \longrightarrow
\]
in the homotopy category of complexes of \Rlm s, $\KModR$.  Here $C$
is bounded to the left because both $M$ and $P^{\ast}$ are bounded to
the left.  This induces a distinguished triangle
\[
  \begin{diagram}[labelstyle=\scriptstyle]
    \Hom_{\R}(C,Q)
    & \rTo &
    \Hom_{\R}(P^{\ast},Q)
    & \rTo^{\Hom_{\R}(\pi^{\ast}\mu,Q)} &
    \Hom_{\R}(M,Q)
    & \rTo & \\
  \end{diagram}
\]
which shows that the chain map in the lemma is a quasi-isomorphism if
and only if the complex $\Hom_{\R}(C,Q)$ is exact.

Now, if the complex $Q$ is just a single projective module placed in
degree zero, then the lemma follows from lemma \ref{lem:Q1}.  So in
this case, $\Hom_{\R}(C,Q)$ must be exact.

Hence $C$ is a complex bounded to the left for which the complex
$\Hom_{\R}(C,Q)$ is exact when $Q$ is a single projective module
placed in degree zero.  But then it is classical that $\Hom_{\R}(C,Q)$
is exact when $Q$ is any complex of projective modules.  Indeed, this
follows from an argument analogous to the one which shows that if $X$
is a complex bounded to the left which is exact and $I$ is any complex
of injective modules, then $\Hom_{\R}(X,I)$ is exact.
\end{proof}

As indicated in the introduction, the category of projective \Rlm s is
denoted $\Proj(\R)$, and the corresponding homotopy category of
complexes is denoted $\KProjR$.  So $\KProjR$ has as objects all
complexes of projective \Rlm s, and as morphisms it has homotopy
classes of chain maps.

\begin{Lemma}
\label{lem:Q3}
For each \fp\ \Rlm\ $M$, there is a natural equivalence
\[
  \Hom_{\KProjRsmall}(P^{\ast},-) 
  \simeq \H^0\!\Hom_{\R}(M,-)
\]
of functors on $\KProjR$.
\end{Lemma}

\begin{proof}
I have
\[
  \Hom_{\KProjRsmall}(P^{\ast},-) 
  \simeq \H^0\!\Hom_{\R}(P^{\ast},-)
  \simeq \H^0\!\Hom_{\R}(M,-)
\]
as functors on $\KProjR$, where the first $\simeq$ is classical and
the second $\simeq$ is by lemma \ref{lem:Q2}.
\end{proof}

\begin{Proposition}
\label{pro:compact}
For each \fp\ \Rlm\ $M$, the complex $P^{\ast}$ from construction
\ref{con:P} is a compact object of $\KProjR$.
\end{Proposition}

\begin{proof}
This is clear from lemma \ref{lem:Q3}, since the functor 
$\H^0\!\Hom_{\R}(M,-)$ respects set indexed coproducts because $M$ is
\fp. 
\end{proof}

\section{Compact generators}
\label{sec:generators}

\begin{Setup}
\label{set:generators}
In this section, $\R$ is a coherent ring (that is, it is both left-
and right-coherent) for which each flat \Rlm\ has finite projective
dimension.
\end{Setup}

\begin{Remark}
\label{rmk:uniform}
Note that there is an integer $N$ so that the projective dimension of
each flat \Rlm\ $F$ satisfies $\pd F \leq N$.  For otherwise, if
there were flat \Rlm s of arbitrarily high, finite projective
dimension, then the coproduct of such modules would be a flat module
of infinite projective dimension.
\end{Remark}

\begin{Construction}
\label{con:generatingset}
For each \fp\ \Rlm\ $M$, take the complex $P^{\ast}$ from construction
\ref{con:P}, and consider the collection of all suspensions $\Sigma^i
P^{\ast}$.

There is only a set (as opposed to a class) of isomorphism classes of
such modules $M$, so there is also only a set of isomorphism classes
in $\KProjR$ of complexes of the form $\Sigma^i P^{\ast}$.  Let the
set $\generatingset$ consist of one object from each such isomorphism
class. 
\end{Construction}

\begin{Theorem}
\label{thm:generators}
The category $\KProjR$ is a \cg\ \tr\ category with $\generatingset$
as a set of compact generators.
\end{Theorem}

\begin{proof}
Each complex $P^{\ast}$ is a compact object of $\KProjR$ by
proposition \ref{pro:compact}, so the same holds for each complex
$\Sigma^i P^{\ast}$ in $\generatingset$.  It remains to show that
$\generatingset$ is a set of generators.  So suppose that $Q$ in
$\KProjR$ has $\Hom_{\KProjRsmall}(G,Q) = 0$ for each $G$ in
$\generatingset$.  I must show $Q \cong 0$ in $\KProjR$.

First, I can consider construction \ref{con:P} with $M$ equal to $\R$,
viewed as an \Rlm.  The corresponding complex $P^{\ast}$ has
suspensions $\Sigma^i P^{\ast}$, and by the construction of
$\generatingset$ each $\Sigma^i P^{\ast}$ is isomorphic to a complex
in $\generatingset$, so  $\Hom_{\KProjRsmall}(\Sigma^i P^{\ast},Q)$
is zero.  Hence
\begin{eqnarray*}
  0 & =     & \Hom_{\KProjRsmall}(\Sigma^i P^{\ast},Q) \\
    & \cong & \Hom_{\KProjRsmall}(P^{\ast},\Sigma^{-i}Q) \\
    & \cong & \H^0\!\Hom_{\R}(A,\Sigma^{-i}Q) \\
    & \cong & \H^{-i}Q,
\end{eqnarray*}
where the second $\cong$ is by lemma \ref{lem:Q3}.  So $Q$ is exact.

Secondly, let me show that for each $j$, the $j$'th cycle module
$\Z^j\!Q$ of $Q$ is flat.  It is clearly enough to do this for
$\Z^0\!Q$.  I shall use the criterion of \cite[chp.\ VI, exer.\
6]{CE}.  So suppose that $a_1, \ldots, a_m$ in $\R$ and $z_1, \ldots,
z_m$ in $\Z^0\!Q$ satisfy the relation
\begin{equation}
\label{equ:a_z_relation}
  \sum_s a_s z_s = 0.
\end{equation}

Consider the \fg\ submodule $M = \R z_1 + \cdots + \R z_m$ of
$\Z^0\!Q$.  Since $\Z^0\!Q$ is a submodule of $Q^0$, so is $M$, and as
$M$ is \fg\ while $Q^0$ is projective and $\R$ coherent, it follows
that $M$ is finitely presented.  So $M$ is among the modules
considered in construction \ref{con:P}, and there is a corresponding
complex $P^{\ast}$.  As above, by the construction of
$\generatingset$ the complex $P^{\ast}$ is isomorphic to a complex in
$\generatingset$, so $\Hom_{\KProjRsmall}(P^{\ast},Q)$ is zero.  Hence
\[
  0 = \Hom_{\KProjRsmall}(P^{\ast},Q)
  \cong \H^0\!\Hom_{\R}(M,Q)
\]
by lemma \ref{lem:Q2}.

So each homomorphism $M \longrightarrow Q^0$ for which the composition
$M \longrightarrow Q^0 \longrightarrow Q^1$ is zero factors through
$Q^{-1} \longrightarrow Q^0$.  In other words, each homomorphism $M
\longrightarrow \Z^0\!Q$ factors through the canonical surjection
$Q^{-1} \stackrel{\sigma}{\longrightarrow} \Z^0\!Q$.  But $M$ is a
submodule of $\Z^0\!Q$, so in particular the inclusion $M
\hookrightarrow \Z^0\!Q$ factors,
\[
  \begin{diagram}[labelstyle=\scriptstyle]
           &               & M                   \\
           & \SW^{f}       & \dInto                \\
    Q^{-1} & \rTo_{\sigma} & \Z^0\!Q \lefteqn{.} \\
  \end{diagram}
\]

Applying $f$ to $\sum_s a_s z_s = 0$ gives $\sum_s a_s f(z_s) = 0$ in
$Q^{-1}$.  But $Q^{-1}$ is projective, hence flat, and so by
\cite[chp.\ VI, exer.\ 6]{CE} there exist $a_{11}, \ldots, a_{mn}$ in
$\R$ and $q_1, \ldots, q_n$ in $Q^{-1}$ so that
\begin{equation}
\label{equ:a_q_relation}
  f(z_s) = \sum_t a_{st}q_t
\end{equation}
and
\begin{equation}
\label{equ:a_a_relation}
  \sum_s a_s a_{st} = 0.
\end{equation}
Applying $\sigma$ to equation \eqref{equ:a_q_relation} gives
\begin{equation}
\label{equ:a_sigma_relation}
  z_s = \sum_t a_{st} \sigma(q_t).
\end{equation}

However, when equation \eqref{equ:a_z_relation} implies the existence
of $a_{11}, \ldots, a_{mn}$ in $\R$ and $\sigma(q_1), \ldots,
\sigma(q_n)$ in $\Z^0\!Q$ so that equations \eqref{equ:a_a_relation}
and \eqref{equ:a_sigma_relation} are satisfied, then \cite[chp.\ VI,
exer.\ 6]{CE} says that $\Z^0\!Q$ is flat as desired.

Finally, note that by remark \ref{rmk:uniform} there is an integer $N$
so that each flat \Rlm\ $F$ has $\pd F \leq N$.  Hence $\pd
\Z^{j+N}\!Q \leq N$ for each $j$.  But there is an exact sequence
\[
  0 \rightarrow \Z^j\!Q \rightarrow Q^j \rightarrow \cdots \rightarrow
  Q^{j+N-1} \rightarrow \Z^{j+N}\!Q \rightarrow 0,
\]
and since $Q^j, \ldots, Q^{j+N-1}$ are projective there follows
$\pd \Z^j\!Q \leq 0$, that is, $\Z^j\!Q$ is projective for each $j$.

So $Q$ is an exact complex of projectives where each cycle module is
also projective.  Hence $Q$ is split exact, and so in particular null
homotopic, so $Q \cong 0$ in $\KProjR$ as desired.
\end{proof}

\section{The subcategory of compact objects}
\label{sec:subcategory}

\begin{Setup}
\label{set:subcategory}
In this section, $\R$ is again a coherent ring for which each flat
$\R$-left-module has finite projective dimension.
\end{Setup}

The \cg\ \tr\ category $\KProjR$ has the full subcategory
$\KProjR^{\c}$ of compact objects.  And the derived category
$\D(\R^{\opp})$ of \Rrm s has the full subcategory $\Df(\R^{\opp})$ of
complexes with bounded cohomology consisting of \fp\ modules.

\begin{Theorem}
\label{thm:subcategory}
There is an equivalence of \tr\ categories
\[
  \KProjR^{\c} \stackrel{\simeq}{\longrightarrow} \Df(\R^{\opp})^{\opp}.
\]
\end{Theorem}

\begin{proof}
Consider again the set $\generatingset$ from construction
\ref{con:generatingset}.  Theorem \ref{thm:generators} says that
$\generatingset$ is a set of compact generators for $\KProjR$.

Let $\sC$ be the full subcategory of $\KProjR$ consisting of objects
which are finitely built from objects $\generatingobj$ in
$\generatingset$.  Let $\sD$ be the full subcategory of
$\KProjRopp$ consisting of objects which are finitely built from
objects of the form $\generatingobj^{\ast}$ with $\generatingobj$
in $\generatingset$.

Each object $\generatingobj$ in $\generatingset$ is a complex of
finitely generated projective modules, so the canonical chain maps
$\generatingobj \rightarrow \generatingobj^{\ast \ast}$ and
$\generatingobj^{\ast} \rightarrow \generatingobj^{\ast \ast \ast}$
are isomorphisms.  Hence
\begin{equation}
\label{equ:equivalences}
  \begin{diagram}[labelstyle=\scriptstyle]
    \sC & \pile{ \rTo^{(-)^{\ast}} \\ \lTo_{(-)^{\ast}} } & \sD^{\opp} \\
  \end{diagram}
\end{equation}
are quasi-inverse equivalences of \tr\ categories.  Indeed, let me
show that this gives the equivalence stated in the theorem: First, the
category $\sC$ consists of the objects finitely built from a set of
compact generators of the \cg\ \tr\ category $\KProjR$, so $\sC$ is
equal to $\KProjR^{\c}$ by the Thomason Localization Theorem,
\cite[thm.\ 2.1]{NeemanDuality}.

Secondly, let me consider the category $\sD$.  It consists of the
objects finitely built from objects of the form
$\generatingobj^{\ast}$ with $\generatingobj$ in $\generatingset$.  By
the definition of $\generatingset$, there is one object $G$ in each
isomorphism class of objects of the form $\Sigma^i P^{\ast}$ with
$P^{\ast}$ coming from construction \ref{con:P}.  So up to
isomorphism, there is one object $G^{\ast}$ in each isomorphism class
of objects of the form $\Sigma^j P$ with $P$ coming from construction
\ref{con:P}.  Recall from construction \ref{con:P} that $P$ is a
projective resolution of the \Rrm\ $M^{\ast}$ which comes from the
\fp\ \Rlm\ $M$.  It follows that $\sD$ consists of the objects
finitely built from projective resolutions of the form $P$.

Now, if $\sD$ had consisted of the objects finitely built from
projective resolutions of {\em all} \fp\ \Rrm s, then $\sD$ would have
been the subcategory of $\KProjRopp$ consisting of projective
resolutions of all complexes with bounded \fp\ cohomology, and it is
classical that this subcategory is equivalent to $\Df(\R^{\opp})$.  So
I would have been done: Equation \eqref{equ:equivalences} would have
given the equivalence stated in the theorem.

As it is, $\sD$ only consists of objects finitely built from
projective resolutions $P$ of \Rrm s of the form $M^{\ast}$
with $M$ a \fp\ \Rlm.  However, this makes no difference because it
turns out that I can finitely build the projective resolution of any
\fp\ \Rrm\ from projective resolutions of the form $P$.

To see this, suppose that $N$ is a \fp\ \Rrm, and let 
\[
  Q = \cdots \rightarrow Q^{-2} \rightarrow Q^{-1}
  \rightarrow Q^0 \rightarrow 0 \rightarrow \cdots
\]
be a projective resolution of $N$.  Since all projective resolutions
of $N$ are isomorphic in $\KProjRopp$, I can suppose that $Q$ consists
of \fg\ projective \Rrm s.  

Now
\[
  \widetilde{Q} = \cdots \rightarrow Q^{-4} \rightarrow Q^{-3}
  \rightarrow Q^{-2} \rightarrow 0 \rightarrow \cdots
\]
is the double suspension of a projective resolution of $\Z^{-1}\!Q$,
the $(-1)$'st cycle module of $Q$, and the complex $Q$ is finitely
built from $Q^0$ and $Q^{-1}$ (viewed as complexes concentrated in
degree zero) along with $\widetilde{Q}$.

Both $Q^0$ and $Q^{-1}$ are projective resolutions of the form
$P$, since they are both projective resolutions of modules of the form
$M^{\ast}$, namely, they are resolutions of $(Q^{0 \ast})^{\ast} \cong
Q^0$ and $(Q^{-1 \ast})^{\ast} \cong Q^{-1}$.  

And $\widetilde{Q}$ is the double suspension of a projective
resolution of the form $P$ because $\Z^{-1}\!Q$ has the form
$M^{\ast}$ for a \fp\ \Rlm\ $M$.  To see this, complete $Q^{0
\ast} \rightarrow Q^{-1 \ast}$ with its cokernel, 
\[
  Q^{0 \ast} \rightarrow Q^{-1 \ast} \rightarrow M \rightarrow 0.
\]
Here $M$ is \fp\, and $M^{\ast}$ sits in the exact sequence
\[
  0 \rightarrow M^{\ast} \rightarrow Q^{-1 \ast \ast} 
  \rightarrow Q^{0 \ast \ast}.
\]
But $Q^0$ and $Q^{-1}$ are \fg, so up to isomorphism the last map here
is just $Q^{-1} \rightarrow Q^0$, so up to isomorphism, the kernel
$M^{\ast}$ is just the kernel of $Q^{-1} \rightarrow Q^0$, that is, it
is $\Z^{-1}\!Q$.  So $\Z^{-1}\!Q$ has the form $M^{\ast}$.
\end{proof}

\section{The dualizing complex case}
\label{sec:dualizing}

\begin{Setup}
\label{set:dualizing}
In this section, $k$ is a field, $\R$ is a $k$-algebra which is
left-coherent and right-noetherian, $\B$ is a left-noetherian
$k$-algebra, and ${}_{\B}D_{\R}$ is a dualizing complex over $\B$ and
$\R$.
\end{Setup}

See \cite[def.\ 1.1]{YekutieliZhang} for the definition of dualizing
complexes. 

\begin{Theorem}
\label{thm:dualizing}
There is an equivalence of \tr\ categories
\[
  \KProjR^{\c} \stackrel{\simeq}{\longrightarrow} \Df(\B).
\]
\end{Theorem}

\begin{proof}
Since there is a dualizing complex ${}_{\B}D_{\R}$ between $\B$ and
$\R$, each flat \Rlm\ has finite projective dimension by
\cite{PJfdpd}.  Moreover, $\R$ is clearly coherent.  So
section \ref{sec:subcategory} applies to $\R$, and theorem
\ref{thm:subcategory} gives an equivalence
\[
  \KProjR^{\c} \stackrel{\simeq}{\longrightarrow} \Df(\R^{\opp})^{\opp}.
\]
But existence of ${}_{\B}D_{\R}$ gives an equivalence
\[
  \Df(\R^{\opp})^{\opp} \stackrel{\simeq}{\longrightarrow} \Df(B)
\]
by \cite[prop.\ 1.3(2)]{YekutieliZhang}, and composing the two
equivalences proves the theorem.
\end{proof}

\bigskip

\noindent
{\bf Acknowledgement.}  
I thank Henning Krause for conversations which strongly inspired this
paper, and for providing me with the crucial trick in the proof of
theorem \ref{thm:subcategory}.

The diagrams were typeset with Paul Taylor's {\tt diagrams.tex}.

\end{document}